\documentclass[12pt]{article}

\usepackage{amssymb}
\usepackage{amsmath}

\begin{document}

% Use the option doublespacing or reviewcopy to obtain double line spacing

%\documentclass[doublespacing]{elsart}

\title{Densit\'e de points et minoration de hauteur}

\author{Nicolas Ratazzi}
\date{}
\maketitle

\vspace{.5cm}

\begin{center} \textit{Universit\'e Paris 6, Projet th\'eorie des nombres, UMR 7586, case 247, 4 place Jussieu, Institut de math\'ematiques, 75252 Paris, FRANCE}
\end{center}

\hrule

\vspace{.5cm}

\noindent \textbf{Abstract}

\vspace{.3cm}

\noindent We obtain a lower bound for the normalised height of a non-torsion subvariety $V$ of a C.M. abelian variety. This lower bound is optimal in terms of the geometric degree of $V$, up to a power of a ``log''. We thus extend the results of F. Amoroso and S. David on the same problem on a multiplicative group $\mathbb{G}_m^n$. We prove furthermore that the optimal lower bound (conjectured by S. David and P. Philippon) is a corollary of the conjecture of S. David and M. Hindry on the abelian Lehmer's problem. We deduce these results from a density theorem on the non-torsion points of $V$.

\vspace{.5cm}

\hrule

\vspace{.5cm}

\noindent \textit{keyword :} abelian varieties, normalised height, Lehmer problem


\footnotetext[1]{ \textit{Email adress : } ratazzi@math.jussieu.fr (Nicolas Ratazzi)}

\vspace{.5cm}

\hrule

\newcounter{ndefinition} 
\newcommand{\defi}{\addtocounter{ndefinition}{1}{\noindent \textbf{D{\'e}finition \thendefinition\ }}}
\newcounter{nrem}
\newcommand{\rem}{\addtocounter{nrem}{1}{\noindent \textbf{Remarque \thenrem\ }}}
\newcommand{\demo}{\noindent \textbf{D\'emonstration} }

\newtheorem{conj}{Conjecture} 
\newtheorem{lemme}{Lemme} 
\newtheorem{propo}{Proposition} 
\newtheorem{coro}{Corollaire} 
\newtheorem{theo}{Th{\'e}or{\`e}me}

\section{Introduction}
\noindent On sait depuis les travaux de Philippon \cite{phi1} \cite{phi2} \cite{phi3}, puis Bost, Gillet, Soul\'e \cite{BGS} dans le cadre de l'intersection arithm\'etique, comment d\'efinir la hauteur des vari\'et\'es projectives ; l'id\'ee \'etant de consid\'erer un point comme une vari\'et\'e de dimension z\'ero et de g\'en\'eraliser ceci en dimension sup\'erieure. De m\^eme que dans le cas des points, on sait pour les vari\'et\'es ab\'eliennes munies d'un fibr\'e en droites ample et sym\'etrique d\'efinir une hauteur particulierement agr\'eable~: la hauteur canonique $\widehat{h}_L$, ou hauteur normalis\'ee. En dimension z\'ero, il existe un th\'eor\`eme caract\'erisant les points de hauteur normalis\'ee nulle ; c'est un r\'esultat de Kronecker dans le cas de $\mathbb{G}_m$. Philippon \cite{phi3} (dans le cas d'un produit de courbes elliptiques) puis Zhang \cite{zhang2} et David-Philippon \cite{daviphi} dans le cas g\'en\'eral ont montr\'e comment g\'en\'eraliser ce r\'esultat pour caract\'eriser les sous-vari\'et\'es de hauteur normalis\'ee nulle : ce sont les translat\'ees d'un sous-groupe alg\'ebrique par un point de torsion. On dit qu'une telle sous-vari\'et\'e est une sous-vari\'et\'e de torsion. La r\'eponse \`a cette question r\'epond \`a une conjecture de Bogomolov. Ceci \'etant, on peut se demander comment minorer la hauteur normalis\'ee d'une sous-vari\'et\'e de hauteur non-nulle d'une vari\'et\'e ab\'elienne. Dans leur article \cite{daviphi}, David et Philippon ont formul\'e un probl\`eme g\'en\'eral (le probl\`eme 1.7) contenant cette question. On peut notamment faire ressortir de la discussion suivant la formulation de leur probl\`eme l'\'enonc\'e suivant :

\vspace{.3cm}

\begin{conj}\label{conj1} (David-Philippon) Soit $A$ une vari\'et\'e ab\'elienne d\'efinie sur un corps de nombres $k$, munie d'un fibr\'e ample et sym\'etrique $L$. Soit $V$ une sous-vari\'et\'e stricte de $A$ sur $k$, $k$-irr\'eductible et qui n'est pas r\'eunion de sous-vari\'et\'es de torsion, alors, on a l'in\'egalit\'e
\[\frac{\widehat{h}_{L}(V)}{\deg_L(V)}\geq c(A,L)\deg_L(V)^{-\frac{1}{s-\textnormal{dim} V}},\]
\noindent o\`u $s$ est la dimension du plus petit sous-groupe alg\'ebrique contenant $V$, et o\`u $c(A,L)$ est une constante ne d\'ependant que de $A$ et de $L$.
\end{conj}

\vspace{.3cm}

\noindent Dans ce qui suit, on reprend le r\'esultat principal (ainsi que le sch\'ema de d\'emonstration) de Amoroso-David \cite{amodavram} concernant le groupe multiplicatif $\mathbb{G}_m^n$, pour obtenir un r\'esultat analogue dans le cadre des vari\'et\'es ab\'eliennes. En utilisant les r\'esultats de David-Hindry \cite{davidhindry} concernant le probl\`eme de Lehmer ab\'elien, on obtient en corollaire un r\'esultat en direction de la conjecture \ref{conj1}. Ce dernier ne concerne que les vari\'et\'es ab\'eliennes de type C.M., par contre il est essentiellement optimal (\`a un facteur $\log$ pr\`es) en le degr\'e de $V$.

\vspace{.3cm}

\noindent \textbf{Remerciements} Je tiens \`a remercier  M. Hindry pour sa patiente relecture, et je le remercie \'egalement, ainsi que S. David, pour m'avoir encourag\'e \`a \'ecrire cet article. Par ailleurs je souhaite aussi remercier chaleureusement \'E. Gaudron et G. R\'emond pour m'avoir indiqu\'e une erreur dans la preuve du lemme \ref{majoration1} dans une version pr\'eliminaire de cet article.

\subsection{Degr\'e et hauteur}

\noindent Soit $k$ un corps de nombres. On dira que $V$ est une \textit{vari\'et\'e alg\'ebrique} sur $k$ si $V$ est un $k$-sch\'ema de type fini g\'eom\'etriquement r\'eduit. On dira que $G$ est un \textit{groupe alg\'ebrique} sur $k$ si c'est une vari\'et\'e en groupes sur $k$. On dira que $A$ est une \textit{vari\'et\'e ab\'elienne} d\'efinie sur $k$ si c'est un groupe alg\'ebrique connexe propre et lisse sur $k$. Par sous-vari\'et\'e on entendra toujours sous-vari\'et\'e ferm\'ee.

\vspace{.3cm}

\noindent Soient $\mathcal{O}_k$ l'anneau des entiers de $k$, $n$ un entier, et $X$ une vari\'et\'e projective munie d'un plongement $\varphi_L : X \hookrightarrow \mathbb{P}^n_k$ d\'efini par un fibr\'e $L$ tr\`es ample sur $X$. Si $\mathcal{O}(1)$ d\'enote le fibr\'e standard sur $\mathbb{P}^n_{\mathcal{O}_k},$ on a $\varphi_L^*\mathcal{O}(1)_k\simeq L$. On note $\overline{\mathcal{O}(1)}$ le fibr\'e standard muni de la m\'etrique de Fubini-Study. Si $V$ est une sous-vari\'et\'e de $X$, on note $\mathcal{V}_L$ l'adh\'erence sch\'ematique de $\varphi_L(V)$ dans $\mathbb{P}^n_{\mathcal{O}_k}$.

\vspace{.3cm}

\defi On d\'efinit le \textit{degr\'e de la vari\'et\'e $V$} relativement \`a $L$, et on note $\deg_L V$ l'entier 
$\deg_k \left(c_1(\mathcal{O}(1)_k)^{\textnormal{dim}V}\cdot\varphi_L(V)\right)$ o\`u $\deg_k$ est le degr\'e projectif usuel sur $\mathbb{P}^n_k$.

\vspace{.3cm}

\defi On appelle \textit{hauteur de la vari\'et\'e $V$} associ\'ee \`a $L$, et on note $h_L(V)$ 
%le r\'eel $h_{\overline{\mathcal{O}(1)}}(\mathcal{V}_L)$ o\`u $h_{\overline{\mathcal{O}(1)}}(.)$ est 
la hauteur de $\mathcal{V}_L$, au sens de Bost-Gillet-Soul\'e \cite{BGS} p. 945 d\'efinition 3.1.1., associ\'ee au fibr\'e hermitien $\overline{\mathcal{O}(1)}$. Notons que l'on ne normalise pas cette hauteur par le degr\'e $\deg_L(V)$.

\vspace{.3cm}

\rem Par le th\'eor\`eme 3 p. 366 de \cite{soule}, $h_L(V)$ coincide avec la hauteur $h(f_{V,L})$ de Philippon, telle que d\'efinie au paragraphe 2. de \cite{phi3}, o\`u $f_{V,L}$ est une forme \'eliminante de l'id\'eal de d\'efinition de $\varphi_L(V)$ dans $k[X_0,\ldots,X_n]$. (Le terme d'erreur de \cite{soule} disparait du fait du changement de normalisation pour la hauteur de Philippon entre les articles \cite{phi1} et \cite{phi3}).

\vspace{.3cm}

\defi Dans le cas o\`u $X=A$ est une vari\'et\'e ab\'elienne, et o\`u $L$ est en plus sym\'etrique, Philippon \cite{phi3} (dans le cas o\`u $L$ d\'efinit un plongement projectivement normal) puis Zhang \cite{zhang2}, avec des m\'ethodes arakeloviennes, ont montr\'e en utilisant un proc\'ed\'e de limite \`a la N\'eron-Tate, comment d\'efinir une \textit{hauteur canonique}, not\'ee $\widehat{h}_L(.)$, sur l'ensemble des sous-vari\'et\'es de $A$. Cette hauteur v\'erifie notamment : si $X$ est une sous-vari\'et\'e de $A$, de stabilisateur $G_X$, et si $n$ est un entier, alors,
\[\widehat{h}_L\left([n](X)\right)=\frac{n^{2(\textnormal{dim} X+1)}}{\mid \textnormal{ker}\ [n]\cap G_X\mid} \widehat{h}_L(X).\] 

\vspace{.3cm}

\defi Soit $A/k$ une vari\'et\'e ab\'elienne. On dit que $V$ est une \textit{sous-vari\'et\'e de torsion} de $A$ si $V=a+B$ avec $a\in A_{\textnormal{tors}}$ et $B$ une sous-vari\'et\'e ab\'elienne de $A$. On dit que c'est une \textit{sous-vari\'et\'e de torsion stricte} de $A$ si c'est une sous-vari\'et\'e de torsion $a+B$ avec $B$ une sous-vari\'et\'e ab\'elienne stricte de $A$.

\vspace{.3cm}

\noindent D'apr\`es les r\'esultats de Philippon \cite{phi3}, David-Philippon \cite{daviphi} et Zhang \cite{zhang2}, on a, si $V$ est une sous-vari\'et\'e de $A/k$ d\'efinie sur une extension finie $K/k$,
\[\widehat{h}_L(V)=0\textnormal{ si et seulement si $V$ est une sous-vari\'et\'e de torsion.}\]

\vspace{.3cm}

\defi Soient $V$ une sous-vari\'et\'e de $A$ sur $k$, et $\theta$ un nombre r\'eel positif. On pose $V(\theta,L)=\left\{x\in V(\overline{k})\ / \ \widehat{h}_L(x)\leq \theta\right\}.$  On d\'efinit alors le \textit{minimum essentiel} de $V$, et on note $\mu_L^{\textnormal{ess}}(V)$ le r\'eel 
\[\mu_L^{\textnormal{ess}}(V)=\inf \left\{\theta>0 \ / \ \overline{V(\theta,L)}=V\ \right\},\]
\noindent o\`u $\overline{V(\theta,L)}$ est la fermeture de Zariski de $V(\theta,L)$ dans $A$.

\subsection{R\'esultats}

\noindent Soient $k$ un corps de nombres, $A/k$ une vari\'et\'e ab\'elienne de dimension $g$, et $L$ un fibr\'e en droites tr\`es ample sur $A$. On d\'emontre le th\'eor\`eme suivant :

\vspace{.3cm}

\begin{theo}\label{theoprinc} Soient $K/k$ une extension finie, et $V$ une sous-vari\'et\'e alg\'ebri\-que de $A$ sur $K$, $K$-irr\'eductible et contenue dans aucune r\'eunion de sous-vari\'et\'es de torsion strictes de $A$. Alors, pour tout r\'eel $\varepsilon>0$, l'ensemble des points $x\in V(\overline{K})$ d'ordre infini pour toute sous-vari\'et\'e ab\'elienne stricte de $A$, et dont la hauteur de N\'eron-Tate relativement \`a $L$ v\'erifie
\[\widehat{h}_L(x)\leq\frac{\widehat{h}_L(V)}{\deg_L V}+\varepsilon\]
\noindent est Zariski dense dans $V$.
\end{theo}

\vspace{.3cm}

\noindent On peut donner deux corollaires \`a ce th\'eor\`eme. Pour cela, on a besoin d'une d\'efinition.

\vspace{.3cm}

\defi Soient $A$ une vari\'et\'e ab\'elienne d\'efinie sur un corps de nombres $k$, $L$ un fibr\'e en droites ample sym\'etrique, et $x\in A(\overline{k}).$ Suivant \cite{davidhindry} d\'efinition 1.2., on appelle \textit{indice d'obstruction de $x$}, et on note $\delta_L(x)$ la quantit\'e
\[\delta_L(x)=\min\left\{\deg_L X^{\frac{1}{\textnormal{codim }X}}\ \ /\ \ x\in X\right\},\]
\noindent o\`u le minimum est pris sur l'ensemble des sous-vari\'et\'es strictes, X, de $A$ sur $k$,  $k$-irr\'eductibles.

\vspace{.3cm}

\noindent On peut voir le point $x\in A(\overline{k})$ comme une sous-vari\'et\'e de $A$, d\'efinie sur $k$, $k$-irr\'eductible, de dimension $0$ et de degr\'e $[k(x):k]$. Avec cette interpr\'etation, la d\'efinition pr\'ec\'edente admet la g\'en\'eralisation suivante : si $V$ est une sous-vari\'et\'e stricte de $A$ sur $k$, $k$-irr\'eductible, on appelle \textit{indice d'obstruction de $V$}, et on note $\delta_L(V)$ la quantit\'e
\[\delta_L(V)=\min\left\{\deg_L X^{\frac{1}{\textnormal{codim }X}}\ \ /\ \ V\subset X\right\},\]
\noindent o\`u le minimum est pris sur l'ensemble des sous-vari\'et\'es strictes, X, de $A$ sur $k$, $k$-irr\'eductibles.

\vspace{.3cm}

\rem On a par d\'efinition, $1\leq\delta_L(V)\leq \deg_L(V)^{\frac{1}{\textnormal{codim }V}}$. Dans le cas o\`u $V$ est de dimension $0$, on retrouve ainsi le lemme 1.3. de \cite{davidhindry}.

\vspace{.3cm}

\noindent En utilisant le r\'esultat de \cite{davidhindry} concernant le probl\`eme de Lehmer pour les vari\'et\'es ab\'eliennes de type  C.M., on peut alors montrer le r\'esultat suivant :

\vspace{.3cm}

\begin{coro}\label{cor1}  Supposons de plus que $A$ est de type C.M. Soit $V$ une sous-vari\'et\'e alg\'ebrique stricte de $A$ sur $k$, $k$-irr\'eductible et contenue dans aucune r\'eunion de sous-vari\'et\'es de torsion strictes de $A$. Alors, l'ensemble
\[\left\{x\in V(\overline{k})\ \ /\ \ \widehat{h}_L(x)\leq \frac{c(A,L)}{\delta_L(V)}\left(\frac{\log\log(3\delta_L(V))}{\log(3\delta_L(V))}\right)^{\kappa(g)}\ \right\},\]
\noindent n'est pas Zariski dense dans $V$. Ici $c(A,L)$ est une constante ne d\'ependant que de $A$ et de $L$, et $\kappa(g)$ est une constante effectivement calculable ne d\'ependant que de $g$ (par exemple $\kappa(g)=(2g(g+1)!)^{g+2}$ convient).
\end{coro}

\vspace{.3cm}

\rem Ce corollaire est une cons\'equence formelle du th\'eor\`eme \ref{theoprinc} et du th\'eor\`eme principal de \cite{davidhindry}. En particulier, toute am\'elioration dans la direction de la conjecture de Lehmer ab\'elienne, am\'eliore d'autant le corollaire. Le meilleur r\'esultat possible correspondrait au cas o\`u $A/k$ est une vari\'et\'e ab\'elienne quelconque, et o\`u l'on peut prendre $\kappa(g)=0$, i.e., \`a la conjecture de Lehmer ab\'elienne telle que \'enonc\'ee dans \cite{davidhindry}.

\vspace{.3cm}

\noindent Dans ce qui suit, si $A$ est une vari\'et\'e ab\'elienne, on suppose donn\'e avec $A$ une isog\'enie avec un produit de vari\'et\'es ab\'eliennes simples $\prod A_i^{n_i}$. On suppose de plus que la vari\'et\'e produit est munie du fibr\'e associ\'e au plongement 
\[A=\prod_{i=1}^nA_i^{r_i}\hookrightarrow \prod_{i=1}^n\mathbb{P}_{n_i}^{r_i}\overset{\textnormal{Segre}}{\hookrightarrow} \mathbb{P}^N,\]
\noindent les $A_i$ \'etant plong\'ees dans $\mathbb{P}_{n_i}$ par des fibr\'es $L_i$ amples et sym\'etriques. Si $L$ est un fibr\'e en droites sym\'etrique ample sur $A$, on notera par $c(A,L)$ une constante ne faisant intervenir que ces donn\'ees.

\vspace{.3cm}

\begin{coro}\label{corfinal} Si $A$ est de type C.M., $L$ un fibr\'e en droites ample et sym\'etri\-que de $A$, et si $V$ est une sous-vari\'et\'e alg\'ebrique stricte de $A$ sur $k$, $k$-irr\'eductible et qui n'est pas r\'eunion de sous-vari\'et\'es de torsion, alors, on a l'in\'egalit\'e
\[\frac{\widehat{h}_{L}(V)}{\deg_L(V)}\geq \mu^{\textnormal{ess}}_L(V)\geq c(A,L)\deg_L(V)^{-\frac{1}{s-\textnormal{dim} V}}\left(\log(3\deg_L(V))\right)^{-\kappa(s)},\]
\noindent o\`u $s$ est la dimension du plus petit sous-groupe alg\'ebrique contenant $V$.
\end{coro}

\vspace{.3cm}

\rem Soit $P\in A(\overline{k})$. En notant $V$ la sous-vari\'et\'e de $A$ sur $k$ obtenue \`a partir de $P$ en rajoutant tous ses conjugu\'es, on constate que le r\'esultat obtenu est bien une g\'en\'eralisation d'un \'enonc\'e du type Lehmer : il s'agit d'un \'enonc\'e de nature arithm\'etique faisant intervenir le degr\'e d'un corps de d\'efinition de $V$.

\vspace{.3cm}

\rem En fait, il suit des preuves des corollaires \ref{cor1} et \ref{corfinal} le r\'esultat suivant : 
\begin{quotation}
\noindent soient $A/k$ une vari\'et\'e ab\'elienne de dimension $g$, $L$ un fibr\'e en droites sym\'etrique ample sur $A$. On appelle Lehmer$(A/k,L,\gamma)$ la propri\'et\'e sui\-vante~: il existe des constantes $c(A,L)$ ne d\'ependant que de $A$ et de $L$, et $\gamma(g)$ ne d\'ependant que de $g$, telles que pour tout point $x\in A(\overline k)$ qui est d'ordre infini modulo toute sous-vari\'et\'e ab\'elienne stricte de $A$, on a  
\[\widehat{h}_L(x)\geq c(A,L)\delta(x)^{-\gamma(g)}.\]

\vspace{.3cm}

\noindent On note ensuite Minorant$(A/k,L,V,\gamma)$ l'\'enonc\'e : il existe des constantes $c'(A,L)$ ne d\'ependant que de $A$ et de $L$, et $\gamma(g)$ ne d\'ependant que de $g$, telles que si $V$ est une sous-vari\'et\'e alg\'ebrique stricte de $A$ sur $k$, $k$-irr\'eductible et qui n'est pas r\'eunion de sous-vari\'et\'es de torsion, alors, on a l'in\'egalit\'e
\[\frac{\widehat{h}_{L}(V)}{\deg_L(V)}\geq c'(A,L)\deg_L(V)^{-\frac{\gamma(g)}{s-\textnormal{dim} V}},\]
\noindent o\`u $s$ est la dimension du plus petit sous-groupe alg\'ebrique contenant $V$.

\vspace{.3cm}

\noindent Avec ces notations, on a 
\[\textnormal{Lehmer}(A/k,L,\gamma) \Rightarrow \textnormal{ Minorant}(A/k,L,V,\gamma).\]
\end{quotation}

\vspace{.3cm}

\noindent De plus, et avec les notations de la remarque pr\'ec\'edente \thenrem, on trouve dans \cite{davidhindry} la conjecture suivante concernant le probl\`eme de Lehmer ab\'elien :

\vspace{.3cm}

\begin{conj} \label{dh}\textnormal{(David-Hindry)} L'assertion \textnormal{Lehmer}$(A/k,L,1)$ est vraie pour toute vari\'et\'e ab\'elienne $A/k$.
\end{conj}

\vspace{.3cm}

\noindent En utilisant la remarque, on en d\'eduit qu'une bonne minoration de la hauteur sur les points entraine une bonne minoration de la hauteur sur toutes les sous-vari\'et\'es. Plus pr\'ecis\'ement, on a : 

\vspace{.3cm}

\begin{coro} La conjecture \ref{dh} implique la conjecture \ref{conj1}.
\end{coro}

\vspace{.3cm}

\noindent La suite est consacr\'ee \`a une d\'emonstration du th\'eor\`eme \ref{theoprinc} et de ces corollaires. 

\section{La proposition cl\'e}

\begin{propo}\label{propcle}\textnormal{(Amoroso-David)} Soient $n$ un entier et $X$ une sous-vari\'et\'e de $\mathbb{P}^n$ sur $k$, $k$-irr\'eductible. Pour tout $\varepsilon>0$, il existe $\delta_0=\delta_0(\varepsilon,X)>0$ v\'erifiant les propri\'et\'es suivantes :
\begin{quotation}
\noindent soient $\delta$ un entier sup\'erieur \`a $\delta_0$, et $Y$ une sous-vari\'et\'e de $\mathbb{P}^n$ sur $k$ ne contenant pas $X$. Si 
\[\log \deg_k (Y)\leq \frac{\delta\varepsilon}{4\textnormal{dim}X},\]
\noindent alors il existe $x\in (X\setminus Y)(\overline{k})$ tel que 
\[h_{\mathcal{O}(1)}(x)\leq\frac{h_{\mathcal{O}(1)}(X)}{\deg_k (X)}+\varepsilon\ \ \textnormal{ et }\ \ [k(x):k]\leq (\deg_k (X))\delta^{\textnormal{dim}X}.\]
\end{quotation}
\end{propo}
\demo C'est la proposition 2.1. de \cite{amodavram} : cette derni\`ere est simplement \'enonc\'ee avec $Y$ une hypersurface, mais la preuve dans le cas g\'en\'eral reste mot pour mot la m\^eme.  \hfill $\Box$

\vspace{.3cm}

\noindent En utilisant un proc\'ed\'e de limite, on en d\'eduit :

\vspace{.3cm}

\begin{coro}\label{cle}Soit $V$ une sous-vari\'et\'e de $A$ sur $k$, $k$-irr\'eductible. Pour tout $\varepsilon_1>0$, il existe $\delta_1=\delta_1(\varepsilon_1,V,A,L)>0$ v\'erifiant les propri\'et\'es suivantes :
\begin{quotation}
\noindent soient $\delta$ un entier sup\'erieur \`a $\delta_1$, et $W$ une sous-vari\'et\'e de $A$ sur $k$ ne contenant pas $V$. Si 
\[\log \deg_L W\leq \frac{\delta}{4\textnormal{dim}V},\]
\noindent alors il existe $x\in (V\setminus W)(\overline{k})$ tel que 
\[\widehat{h}_L(x)\leq\frac{\widehat{h}_L(V)}{\deg_L V}+\varepsilon_1\ \ \textnormal{ et }\ \ [k(x):k]\leq c_0(\varepsilon_1,A,L)(\deg_L V)\delta^{\textnormal{dim}V},\]
\noindent o\`u $c_0(\varepsilon_1,A,L)$ est une constante strictement positive ne d\'e\-pen\-dant que de $\varepsilon_1$, $A$ et $L$.
\end{quotation}
\end{coro}
\demo Quitte \`a remplacer $L$ par $L^{\otimes 4}$, on suppose que le plongement $\varphi : A \hookrightarrow \mathbb{P}^n$  associ\'e au fibr\'e en droites tr\`es ample sym\'etrique $L$, est projectivement normal. Soit $p$ un nombre premier et $N$ un entier strictement positif. On consid\`ere le plongement projectif $\psi=\psi_{p,L}$ de $A$, compos\'e des plongements suivants
\[
\begin{array}{ccccccc}
A & \hookrightarrow & A^N 			 & \hookrightarrow & (\mathbb{P}^n)^N & \underset{\textnormal{Segre}}{\hookrightarrow} & \mathbb{P}^{(n+1)^N-1}\\
x & \mapsto         & (x,[p]x,\ldots,[p^{N-1}]x) &                 &          &                                          &
\end{array}
\]
\noindent Il s'agit du \textit{plongement enroul\'e} d\'efini dans \cite{phi3} paragraphe 3. En notant $h_{\psi}$ la hauteur associ\'ee \`a ce plongement, la proposition 7. de \cite{phi3} nous dit que 
\[\deg_{\psi}(V)=\left(\frac{p^{2N}-1}{p^2-1}\right)^{\textnormal{dim}V}\deg_L(V),\ \textnormal{ et }\ \  \widehat{h}_{\psi}(V)=\left(\frac{p^{2N}-1}{p^2-1}\right)^{\textnormal{dim}V}\widehat{h}_L(V).\]
\noindent Ainsi, en appliquant la proposition 9. de \cite{phi3}, on en d\'eduit qu'il existe un r\'eel $c_p>0$ ind\'ependant de $N$, tel que 
\begin{equation}\label{h1}
\textnormal{\Huge{$\mid$}} \left(\frac{p^{2N}-1}{p^2-1}\right)\frac{\widehat{h}_L(V)}{\deg_L(V)}-\frac{h_{\mathcal{O}(1)}(\psi(V))}{\deg_k(\psi(V))}\textnormal{\Huge{$\mid$}} \leq 8c_pN.
\end{equation}
\noindent Soit maintenant $\varepsilon_1>0$. On fixe $p=3$ par exemple, et on choisit $N=N(\varepsilon_1,A,L)$ le plus petit entier tel que 
\[\frac{16c_pN+1}{\left(\frac{p^{2N}-1}{p^2-1}\right)}\leq\varepsilon_1.\]
\noindent On va appliquer la proposition pr\'ec\'edente \ref{propcle} avec $\varepsilon=1$, $X=\varphi(V)$ et $Y=\varphi(W)$. On choisit $\delta\geq\delta_1(\varepsilon_1,V,A,L)=\max\left\{8(\textnormal{dim} V)^2N\log p,\ \delta_0(1,X)\right\}$, et on suppose que 
\[\log\deg_L(W)\leq\frac{\delta}{4\textnormal{dim}X}.\]
\noindent Avec ces choix, on a
\begin{align*}
\log\deg_k Y 	& = 	\textnormal{dim} V\log \left(\frac{p^{2N}-1}{p^2-1}\right)+\log\deg_L W\\
		& \leq	2N\textnormal{dim} V\log p + \log\deg_L W\\
		& \leq	2N\textnormal{dim} V\log p + \frac{\delta}{4\textnormal{dim} V} \leq \frac{\delta}{2\textnormal{dim}V}\ \ \textnormal{ par choix de }\delta_1.
\end{align*}
\noindent La proposition \ref{propcle} nous dit qu'il existe $y\in \left(\varphi(V)\setminus\varphi(W)\right)(\overline{k})$ tel que :
\[h_{\mathcal{O}(1)}(y)\leq \frac{h_{\mathcal{O}(1)}(X)}{\deg_k (X)}+1\leq 8c_pN+\frac{p^{2N}-1}{p^2-1}\frac{\widehat{h}_L(V)}{\deg_L (V)}+1,\]
\noindent et tel que
\[[k(y):k]\leq \deg_k(X)\delta^{\textnormal{dim}V}\leq \left(\frac{p^{2N}-1}{p^2-1}\right)^{\textnormal{dim} V}\deg_L(V)\delta^{\textnormal{dim}V}.\]
\noindent Par d\'efinition, il existe $x\in V\setminus W$ tel que 
\[y=\varphi(x), \text{ et tel que }\mid h_{\mathcal{O}(1)}(y)-\frac{p^{2N}-1}{p^2-1}\widehat{h}_L(x)\mid\leq 8c_pN.\]
\noindent On en d\'eduit 
\[\widehat{h}_L(x)\leq \frac{\widehat{h}_L(V)}{\deg_L V}+\frac{16c_pN+1}{\left(\frac{p^{2N}-1}{p^2-1}\right)},\ \textnormal{ et } [k(x):k]\leq \left(\frac{p^{2N}-1}{p^2-1}\right)^{2\textnormal{dim} V}\delta^{\textnormal{dim} V}\deg_L V.\]
\noindent Le choix de $N$ permet de conclure.  \hfill $\Box$

\section{Un lemme de majoration}\label{nouveau}
\noindent Dans ce qui suit, $A/k$ est une vari\'et\'e ab\'elienne d\'efinie sur un corps de nombres $k$, et $L$ est un fibr\'e en droites sym\'etrique tr\`es ample sur $A$. On suppose $k$ plong\'e dans $\mathbb{C}$. On commence par rappeler un r\'esultat classique concernant le corps de d\'efinition d'une sous-vari\'et\'e ab\'elienne de $A$.

\vspace{.3cm}

\begin{lemme}\label{defi} Il existe une extension F/k finie, ne d\'ependant que de $A$, et notamment de degr\'e major\'e par une constante ne d\'ependant que de $A$, telle que toute sous-vari\'et\'e ab\'elienne de $A$ soit d\'efinie sur $F$.
\end{lemme}
\demo Le groupe $\textnormal{End}_{\mathbb{C}}(A)$ est un $\mathbb{Z}$-module de type fini. Il existe donc une extension finie $F/k$ ne d\'ependant que de $A$ telle que $\textnormal{End}_{\mathbb{C}}(A)=\textnormal{End}_{F}(A)$. Soit maintenant $B$ une sous-vari\'et\'e ab\'elienne de $A$. Par le th\'eor\`e\-me d'irr\'eductibilit\'e de Poincar\'e, il existe une sous-vari\'et\'e ab\'elienne $C$ de $A$ et une isog\'enie $\varphi : A\rightarrow B\times C.$ En notant $pr_1$ la projection de $B\times C$ sur $B$, et $i$ l'inclusion de $B$ dans $A$, on a : $f=i\circ pr_1\circ\varphi\in \textnormal{End}_F(A)$. Or $B=f(A)$, donc $B$ est d\'efinie sur $F$.   \hfill $\Box$

\vspace{.3cm}

\rem En fait on peut prendre $F/k$ de degr\'e inf\'erieur \`a $3^{16\textnormal{dim}A^4}$ (cf. \cite{masserwust} lemma 2.2.).

\vspace{.3cm}

\noindent Au vu de ce lemme, on supposera dans toute la suite que toutes les sous-vari\'et\'es ab\'eliennes de $A/k$ sont d\'efinies sur $k$.

\vspace{.3cm}

\begin{lemme}\label{torsion}Soient $d$ un entier, et $x\in A$ un point rationnel sur une extension de degr\'e inf\'erieur \`a $d$. Si $x$ est un point de torsion modulo une sous-vari\'et\'e ab\'elienne stricte $B$ de $A$, alors on peut \'ecrire $x=y+\xi$ avec $y\in B$ et $\xi\in A(F)_{\textnormal{tors}}$ o\`u $F$ est une extension de degr\'e inf\'erieur \`a $c_1(A)d$, $c_1(A)$ \'etant une constante ne d\'ependant que de $A$.
\end{lemme}
\demo On note $\pi : A \rightarrow A/B=C$. On sait (cf. \cite{bertrand}) que l'on peut construire une sous-vari\'et\'e ab\'elienne $C'$ de $A$ telle que $A=B+C'$, et telle que $\textnormal{Card}\left(B\cap C'\right)\leq c_1(A)$ pour une constante $c_1(A)$ ne d\'ependant que de $A$. Notons $\pi'=\pi_{\mid C'}$ l'isog\'enie de $C'$ vers $C$, et posons $K=k(x)$. On peut \'ecrire $x=b+c'$ avec $b\in B$ et $c'\in C'$. On a $\pi(x)=\pi'(c')\in C(K)_{\textnormal{tors}}$. L'application $\pi'$ \'etant une isog\'enie, le point $c'$ est de torsion, et il est rationnel sur une extension de $K$ de degr\'e major\'e par $c_1(A)$.
 \hfill $\Box$

\vspace{.3cm}

\begin{lemme}\label{majoration1}Soient $d$ un entier, et $x\in A$ un point rationnel sur une extension de degr\'e inf\'erieur \`a $d$. Si $x$ est de torsion modulo une sous-vari\'et\'e ab\'elienne stricte $B$ de $A$, alors, il existe une sous-vari\'et\'e ab\'elienne $B_x$ stricte de $A$ et un point de torsion $\xi$ d\'efini sur une extension de degr\'e au plus $c_1(A)d$ de $k$, tels que $x\in\left(B_x +\xi\right),$ et tels que 
\[\deg_L B_x\leq c_2(A,L)d^{c_3(A)}\max\{1,\hat{h}_L(x)\}^{c_4(A)}, \ \ \textnormal{o\`u}\]
$c_2$ ne d\'epend que de $A$ de $L$, et o\`u $c_3$, $c_4$ sont des constantes ne d\'ependant que de $A$.
\end{lemme}
\demo On note $G_x$ le plus petit sous-groupe alg\'ebrique contenant $x$, et on note $B_x=G_x^0$ sa composante neutre. Par hypoth\`ese sur $x$, la sous-vari\'et\'e ab\'elienne $B_x$ est strictement incluse dans $A$, et $x$ est de torsion modulo $B_x$. Ainsi, le lemme \ref{torsion} entraine que $x\in\left(B_x+\xi\right)$ avec $\xi$ point de torsion d\'efini sur une extension de degr\'e au plus $c_1(A)d$ de $k$. Il reste \`a voir que le degr\'e de $B_x$ est major\'e comme on veut. On va pour cela utiliser l'article \cite{bertranduke} de Bertrand. En suivant les notations de cet article, on note $\mathcal{H}(A)$ l'ensemble des classes d'isomorphismes de sous-vari\'et\'es ab\'eliennes de $A$. La proposition 1.(ii) de \cite{bertranduke} nous assure que cet ensemble est fini. De plus, si $K/k$ est une extension de degr\'e $d$ telle que le point $x$ est $K$-rationnel, en utilisant le theorem p. 154 de \cite{masser}, on note que le cardinal de $A(K)_{\textnormal{tors}}$ est major\'e par $c_8'(A)d^{7 \textnormal{dim }A}$. Il en va donc de m\^eme du cardinal de $Y(K)_{\textnormal{tors}}$ pour toute sous-vari\'et\'e ab\'elienne $Y$ de $A$. Par le lemme 2, pour toute sous-vari\'et\'e ab\'elienne $Y$ de $A$, le cardinal du sous-groupe de torsion de $(A/Y)(K)$ divise une quantit\'e major\'ee par $c_9'(A)d^{c_{10}'(A)}$. Avec les notations de \cite{bertranduke} poroposition 1. (i), on peut donc prendre $\nu(A/K)\leq c_9'(A)d^{c_{10}'(A)}$. En appliquant maintenant le Corollary p.239 et la remark 2. (i) p.231 de \cite{bertranduke}, et avec ses notations, on obtient la majoration
\[\deg_L B_x\leq c(A,K)^{-1}c_9''(A,L)d^{c_{11}'(A)}\max\{1,\hat{h}_L(x)\}^g.\]
\noindent Il nous suffit maintenant de montrer que $c(A,K)\geq c_{11}''(A,L)d^{-c_{10}''(A)}$. L'en\-semble $\mathcal{H}(A,K)=\mathcal{H}(A,k)$ de classes de K-isomorphismes de sous-vari\'et\'es ab\'eliennes de $A$ est fini (et ne d\'epend que de $A$ et $k$ par le lemme \ref{defi})~: 
\[\mathcal{H}(A,k)=\left\{B_1,\ldots,B_{n(A)}\right\}.\]
\noindent On note $b_i$ le degr\'e de la polarisation sur $B_i$ d\'eduite de $(A,L)$. On pose $c_{12}'(A,L)=\min b_i>0$. La remarque suivant le corollary p.239 de \cite{bertranduke} nous indique alors (en mettant dans une m\^eme constante $c_2(A,L)$ la d\'ependance en $h_{\textnormal{Falt}}(A)$ et en $c_{12}'(A,L))$ que 
\[c(A,K)\geq c_{11}''(A,L)d^{-c_{10}''(A)}.\]
 \hfill $\Box$

\vspace{.3cm}

\rem La preuve de ce lemme nous permet m\^eme de sp\'ecifier $B_x$~: on peut prendre pour $B_x$ la composante neutre du plus petit sous-groupe alg\'ebrique contenant le point $x$.

\vspace{.3cm}

\rem Dans son article \cite{gael}, R\'emond obtient une version plus fine de ce lemme \ref{majoration1}.

\vspace{.3cm}

\noindent Soit $D$ un entier. On d\'efinit une sous-vari\'et\'e de $A$ sur $k$, not\'ee $Y(D,d)$, par  : 
\[Y(D,d)=\bigcup_{B}\left(\bigcup_{\xi} \left(B+\xi\right)\right),\]
o\`u $B$ d\'ecrit l'ensemble $\mathcal{E}_L(D)$ des sous-vari\'et\'es ab\'eliennes strictes de $A$ de degr\'e (relativement \`a $L$) inf\'erieur \`a $D$, et $\xi$ d\'ecrit l'ensemble des points de torsion de $A$ d\'efinis sur une extension de degr\'e au plus $c_1(A)d$ de $k$.

\vspace{.3cm}

\begin{lemme}\label{card} Il existe une constante $c_5(A,L)$ telle que 
\[\textnormal{Card } \mathcal{E}_L(D)\leq c_5(A,L) D^{(2\textnormal{dim}A)^2}.\]
\end{lemme}
\demo Le fibr\'e $L$ est tr\`es ample, donc il d\'efinit une forme de Riemann $H_L$ sur $t_{A(\mathbb{C})}$, telle que la forme symplectique $E_L=\textnormal{Im}H_L$, est \`a valeurs enti\`eres sur le r\'eseau des p\'eriodes $\Omega_{A(\mathbb{C})}$. En utilisant essentiellement le th\'eor\`eme de Riemann-Roch pour les vari\'et\'es ab\'eliennes, la proposition 3 p.269 de \cite{philbert} nous indique que, si $B$ est une sous-vari\'et\'e ab\'elienne de $A$, alors 
\[\deg_L B=(\textnormal{dim}B)!\textnormal{Vol}_{E_L}\Omega_{B(\mathbb{C})},\]
\noindent o\`u le volume est relatif \`a la norme $\mid\mid\cdot\mid\mid$ induite par la forme bilin\'eaire sym\'etrique d\'efinie positive $e_L$ donn\'ee par $e_L(x,y)=E_L(ix,y)$. En notant $c'_6(A,L)$ le volume de la boule unit\'ee de $\mathbb{R}^{2\textnormal{dim}A}$ pour cette norme, le th\'eor\`eme des minimas successifs de Minkowski nous permet d'en d\'eduire qu'il existe une base de $\Omega_{B(\mathbb{C})}$ form\'ee d'\'el\'ements (not\'es $\{\omega_1,\ldots,\omega_{2\textnormal{dim}B}\}$) appartenant \`a $\Omega_{A(\mathbb{C})}$, tels que 
\[\prod_{i=1}^{2\textnormal{dim}B}\mid\mid\omega_i\mid\mid\leq c'_6(A,L)\deg_L B.\]
\noindent On pose $c_{\min}(A,L)=\min\left\{1,\mid\mid\lambda_i\mid\mid\ /\ \lambda_i\in \Omega_{A(\mathbb{C})}-\{0\}\right\}$, et on note $\omega_{\max}$ le $\omega_i$ de plus grande norme. On a 
\[\mid\mid \omega_{\max}\mid\mid c_{\min}(A,L)^{2g-1}\leq \prod_{i=1}^{2\textnormal{dim}B}\mid\mid\omega_i\mid\mid\leq c'_6(A,L)D.\]
\noindent Ainsi, on tire $\mid\mid\omega_{\max}\mid\mid\leq c'_7(A,L)D,$ ce qui entraine 
\[\textnormal{Card }\mathcal{E}_L(D)\leq c_5(A,L)D^{(2\textnormal{dim}A)^2}.\] \hfill $\Box$

\vspace{.3cm}

\begin{lemme}\label{majoration2} Le degr\'e relativement \`a $L$ de $Y(D,d)$ est major\'e par une expression de la forme $c_6(A,L)D^{c_7(A)}d^{c_8(A)}.$
\end{lemme}
\demo On sait par le theorem p. 154  de \cite{masser} que tout point de torsion de $A(\overline{k})$  d\'efini sur une extension de degr\'e inf\'erieur \`a $c_1(A)d$ est d'ordre major\'e par $c'_8(A)d^{7\textnormal{dim}(A)}.$ On en d\'eduit donc que  l'ensemble des points de torsion d\'efinis sur une extension de degr\'e inf\'erieur \`a $c_1(A)d$ est de cardinal major\'e par

\vspace{.3cm}

\begin{equation}
c'_9(A)d^{(7\textnormal{dim}(A))(2\textnormal{dim}(A)+1)}.\label{eq1}
\end{equation}
\noindent Par ailleurs, on sait majorer le cardinal de $\mathcal{E}_L(D)$ par le lemme pr\'ec\'edent, donc l'additivit\'e du degr\'e nous donne
\begin{equation}\label{eq2}
\deg_L Y(D,d)\leq D\cdot\left(c_5(A,L) D^{(2\textnormal{dim}A)^2}\right)\cdot\left(c'_9(A)d^{(7\textnormal{dim}(A))(2\textnormal{dim}(A)+1)}\right).
\end{equation}
\noindent L'in\'egalit\'e (\ref{eq2}) est bien de la forme voulue.   \hfill $\Box$

\section{Preuve du th\'eor\`eme \ref{theoprinc}}

\noindent Quitte \`a remplacer $V$ par la r\'eunion des $\sigma(V)$ avec $\sigma\in\textnormal{Gal}(\overline{k}/k)$, on peut supposer que $V$ est d\'efinie sur $k$ et $k$-irr\'eductible. Soit $\varepsilon>0$ un r\'eel. On suppose par l'absurde qu'il existe une hypersurface $Z$ de $A$ sur $k$, ne contenant pas $V$ mais contenant tous les points de $V$ d'ordre infini modulo toute sous-vari\'et\'e ab\'elienne stricte de $A$, et tels que la hauteur v\'erifie $\widehat{h}_L(x)\leq \frac{\widehat{h}_L(V)}{\deg_L V}+\varepsilon.$ Soit alors, $\delta$ un entier non nul. On consid\`ere la sous-vari\'et\'e $Y_{\delta}$ de $A$ sur $k$ de codimension sup\'erieure \`a $1$, d\'efinie par
\[Y_{\delta}=Y(D,d),\]
\noindent o\`u $Y(D,d)$ est d\'efinie comme au paragraphe pr\'ec\'edent, et o\`u 
\[ D=c_2(A,L)\left(\frac{\widehat{h}_L(V)}{\deg_L(V)}+\varepsilon\right)^{c_4(A)}\left(c_0(\varepsilon)\deg_L(V)\delta^{\textnormal{dim}V}\right)^{c_3(A)},\]
\noindent et,
\[ d=c_0(\varepsilon)\deg_L(V)\delta^{\textnormal{dim}V}.\]
\noindent Par hypoth\`ese, $V$ n'est contenue dans aucune r\'eunion de sous-vari\'et\'es de torsion strictes de $A$, donc $V\nsubseteq Y_{\delta}$. Il existe $\delta_1(\varepsilon,V,A,L)$ tel que pour tout $\delta\geq\delta_1(\varepsilon,V,A,L)$, on a l'in\'egalit\'e
\[\log \deg_L(Y_{\delta}\cup Z)\leq \frac{\delta}{4\textnormal{dim}V}.\]
\noindent En effet, par le lemme \ref{majoration2}, on sait que 
\[\deg_L Y_{\delta}\leq c_6(A,L)D^{c_7(A)}d^{c_8(A)}.\]
\noindent En rempla\c{c}ant $D$ et $d$ par leurs valeurs, et par additivit\'e du degr\'e, on en d\'eduit l'in\'egalit\'e pour tout $\delta\geq\delta_1(\varepsilon,V,A,L)$ assez grand. On se fixe d\'esormais un tel $\delta$, et on applique le corollaire \ref{cle} \`a $Y_{\delta}\cup Z$. On obtient ainsi un $x\in V\setminus(Y_{\delta}\cup Z)$ tel que 
\[\widehat{h}_L(x)\leq\frac{\widehat{h}_L(V)}{\deg_L V}+\varepsilon\ \ \textnormal{ et }\ \ [k(x):k]\leq c_0(\varepsilon)(\deg_L V)\delta^{\textnormal{dim}V}.\]
\noindent Si $x$ est un point de torsion modulo une sous-vari\'et\'e ab\'elienne stricte de $A$, alors, le lemme \ref{majoration1} et le choix de $D$ dans $Y_{\delta}$ entraine que $x\in Y_{\delta}$. Ceci est impossible, donc par d\'efinition de $Z$, $x$ appartient \`a $Z$. Mais ceci est \'egalement impossible. Ceci conclut par l'absurde. \hfill $\Box$ 

\section{Preuve du corollaire \ref{cor1}}
\noindent On commence par rappeler le th\'eor\`eme de Zhang sur les minimas successifs (cf. \cite{zhang1} theorem 5.2, et \cite{zhang2} theorem 1.10). Plus exactement, on en donne une version affaiblie qui nous suffira, ne faisant intervenir que le minimum essentiel.

\vspace{.3cm}

\begin{theo}\label{zh}\textnormal{(Zhang)} Si $V$ est une sous-vari\'et\'e de $A$ sur $K$, alors
\[\frac{\widehat{h}_L(V)}{\left(\textnormal{dim }V+1\right)\deg_L V}\leq \mu_L^{\textnormal{ess}}(V)\leq \frac{\widehat{h}_L(V)}{\deg_L V}.\]
\end{theo}

\vspace{.3cm}

\noindent On peut maintenant passer \`a la preuve du corollaire \ref{cor1} : soit $\varepsilon>0$, le th\'eor\`eme \ref{theoprinc} nous indique que l'ensemble des $x\in V(\overline{k})$ d'ordre infini modulo toute sous-vari\'et\'e ab\'elienne stricte de $A$ et qui sont de hauteur $\widehat{h}_L(x)\leq \frac{\widehat{h}_L(V)}{\deg_L V}+ \varepsilon$ est Zariski dense dans $V$. En utilisant le th\'eor\`eme \ref{zh} de Zhang, on en d\'eduit que les points $x\in V(\overline{k})$ d'ordre infini modulo toute sous-vari\'et\'e ab\'elienne, et de hauteur $\widehat{h}_L(x)\leq (\textnormal{dim }V+1)\mu_L^{\textnormal{ess}}(V)+\varepsilon$ est Zariski dense dans $V$. En particulier cet ensemble est non-vide. On choisit un \'el\'ement $x$ dedans. En appliquant le th\'eor\`eme 1.5. de \cite{davidhindry} ainsi que la remarque qui suit ce th\'eor\`eme, on en d\'eduit
\[\widehat{h}_L(x)\geq \frac{c(A,L)}{\delta_L(x)}\left(\frac{\log\log(3\delta_L(x))}{\log(3\delta_L(x))}\right)^{\kappa(g)},\]
\noindent o\`u $c(A,L)$ est une constante ne d\'ependant que de $A$ et de $L$, et $\kappa(g)$ est une constante effectivement calculable ne d\'ependant que de $g$ (par exemple $\kappa(g)=(2g(g+1)!)^{g+2}$ convient). Or $\delta_L(x)\leq \delta_L(V)$ car une sous-vari\'et\'e de $A$ contenant $V$ contient $x$. On en d\'eduit 
\[\widehat{h}_L(x)\geq \frac{c(A,L)}{\delta_L(V)}\left(\frac{\log\log(3\delta_L(V))}{\log(3\delta_L(V))}\right)^{\kappa(g)}.\]
Notamment on en conclut
\[\left(\textnormal{dim }V+1\right)\mu_L^{\textnormal{ess}}(V)+\varepsilon\geq \frac{c(A,L)}{\delta_L(V)}\left(\frac{\log\log(3\delta_L(V))}{\log(3\delta_L(V))}\right)^{\kappa(g)}.\]
Ceci termine la preuve en faisant tendre $\varepsilon$ vers $0$. \hfill $\Box$

\section{Preuve du corollaire \ref{corfinal}}

\noindent On commence par prouver le corollaire dans un cas particulier, auquel on se ram\`enera ensuite.

\vspace{.3cm}

\begin{coro}\label{cor3} Si $A=\prod_{i=1}^nA_i^{r_i}$, o\`u les $A_i$ sont simples, est de type C.M., et si $V$ est une sous-vari\'et\'e alg\'ebrique stricte de $A$ sur $k$, $k$-irr\'eductible et qui n'est pas r\'eunion de sous-vari\'et\'es de torsion, alors, on a l'in\'egalit\'e
\[\frac{\widehat{h}_{M}(V)}{\deg_M(V)}\geq \mu^{\textnormal{ess}}_M(V)\geq c(A,M)\deg_M(V)^{-\frac{1}{s-\textnormal{dim} V}}\left(\log(3\deg_M(V))\right)^{-\kappa(s)},\]
\noindent o\`u $s$ est la dimension du plus petit sous-groupe alg\'ebrique contenant $V$, et o\`u $M$ est le fibr\'e en droites ample associ\'e au plongement 
\[A=\prod_{i=1}^nA_i^{r_i}\hookrightarrow \prod_{i=1}^n\mathbb{P}_{n_i}^{r_i}\overset{\textnormal{Segre}}{\hookrightarrow} \mathbb{P}^N,\]
\noindent les $A_i$ \'etant plong\'ees dans $\mathbb{P}_{n_i}$ par des fibr\'es $L_i$ tr\`es amples et sym\'etriques.
\end{coro}
\demo On note $G$ le plus petit sous-groupe alg\'ebrique contenant $V$. On note $G^0$ la composante connexe de l'identit\'e de $G$. C'est une sous-vari\'et\'e ab\'elienne de $A$, et elle est donc isog\`ene \`a $B=\prod_{i=1}^nA_i^{s_i}$ o\`u $0\leq s_i\leq r_i$. On note alors $\pi : A \rightarrow B$ une projection naturelle obtenue par oubli de certaines coordonn\'ees, de sorte que $\pi_{\mid G}$ est une isog\'enie. Montrons maintenant que l'on est dans les conditions d'application du corollaire $\ref{cor1}$ en prenant comme vari\'et\'e ab\'elienne $B$, et comme sous-vari\'et\'e alg\'ebrique $\pi(V)$.

\vspace{.3cm}

Si $\pi(V)$ est inclus dans une r\'eunion de sous-vari\'et\'es de torsion $\bigcup\ (C_i+\xi_j)$ o\`u $\textnormal{dim} \, C_i<\textnormal{dim} B,$ en notant $H$ le plus petit sous-groupe alg\'ebrique contenant $\pi(V)$, on a toujours $\textnormal{dim} H<\textnormal{dim} B$. Ainsi $G_1=G\cap\pi^{-1}(H)$ est un sous-groupe alg\'ebrique strict de $G$ (car $\pi_{\mid G}$ est une isog\'enie), contenant $V$. Ceci est absurde.

\vspace{.3cm}

Si $\pi(V)=B$, alors $V$ est de torsion. Ceci est absurde.

\vspace{.3cm}

\noindent Finalement, $\pi(V)$ est une $k$-sous-vari\'et\'e stricte de $B$, irr\'eductible, et n'est pas incluse dans une r\'eunion de sous-vari\'et\'es de torsion strictes. On peut donc appliquer le corollaire \ref{cor1}. Par ailleurs, la hauteur et le degr\'e sont d\'efinis relativement aux plongements 
\[A=\prod_{i=1}^nA_i^{r_i}\hookrightarrow \prod_{i=1}^n\mathbb{P}_{n_i}^{r_i}\overset{\textnormal{Segre}}{\hookrightarrow} \mathbb{P}^{N_A}, \textnormal{  et }B=\prod_{i=1}^nA_i^{s_i}\hookrightarrow \prod_{i=1}^n\mathbb{P}_{n_i}^{s_i}\overset{\textnormal{Segre}}{\hookrightarrow} \mathbb{P}^{N_B}.\]
\noindent De plus l'application $\overline{\pi} :  \prod_{i=1}^n\mathbb{P}_{n_i}^{r_i}\rightarrow  \prod_{i=1}^n\mathbb{P}_{n_i}^{s_i}$ est la projection lin\'eaire d\'efinie par oubli de coordonn\'ees. Dans ce cas et pour ces plongements on a, 
\[\mu^{\textnormal{ess}}_{M_B}(\pi(V))\leq \mu^{\textnormal{ess}}_M(V),\ \ \textnormal{ et, }\ \ \deg_{M_B}\pi(V)\leq\deg_M(V).\]
\noindent Ceci nous donne
\begin{align*}
\mu^{\textnormal{ess}}_M(V) 	& \geq \mu^{\textnormal{ess}}_{M_B}(\pi(V)),\ \ \ \textnormal{d'o\`u par le corollaire \ref{cor1} et la remarque 2,}\\
			& \geq c(B, M_B) \left(\deg_{M_B}\pi(V)\right)^{-\frac{1}{s-\textnormal{dim} V}}\left(\log\deg_{M_B}\pi(V)\right)^{-\kappa(s)}\\
			& \geq c(B, M_B) \left(\deg_{M}V\right)^{-\frac{1}{s-\textnormal{dim} V}}\left(\log\deg_M(V)\right)^{-\kappa(s)}.\\
			& \geq c'(A, M) \left(\deg_{M}V\right)^{-\frac{1}{s-\textnormal{dim} V}}\left(\log\deg_M(V)\right)^{-\kappa(s)},
\end{align*}

\vspace{.3cm}

\noindent o\`u on a pris pour $c'(A,M)$ le minimum des $c(B,M_B)$ quand $s_i$ varie dans $[\![0,r_i]\!]$. On conclut en appliquant le th\'eor\`eme \ref{zh} de Zhang. \hfill $\Box$
\hfill $\Box$

\vspace{.3cm}

\noindent On donne maintenant la preuve du corollaire \ref{corfinal} : la vari\'et\'e ab\'elienne $A$ est donn\'ee avec une isog\'enie $\rho$ vers $B=\prod_{i=1}^n A_i^{r_i}$. Soit $V$ la sous-vari\'et\'e de $A$ comme dans les hypoth\`eses. On v\'erifie que $W=\rho(V)$ est une sous-vari\'et\'e de $B$ v\'erifiant les m\^emes hypoth\`eses. Il r\'esulte facilement de la preuve de la proposition 14. de \cite{phi3} qu'il existe $c'(A,L)$ tel que
\[\widehat{h}_L(V)\geq c'(A,L) \widehat{h}_M(W).\]
\noindent Ainsi, en appliquant le r\'esultat pr\'ec\'edent, on en d\'eduit presque l'in\'egalit\'e voulue : il faut encore remplacer le degr\'e $\deg_M(W)$ par $\deg_L(V)$. Or 
\[\deg_M(W)=(\deg\rho)\deg_{\rho^*M}(V).\]
\noindent D'autre part $\rho^*M$ et $L$ sont amples, donc on a des in\'egalit\'es
\[c_2(A,L)\deg_{\rho^*M}(V)\geq\deg_LV\geq c_3(A,L)\deg_{\rho^*M}(V).\]
\noindent En injectant ceci dans l'in\'egalit\'e donn\'ee par le corollaire \ref{cor3} pr\'ec\'edent, on peut conclure.

\end{document}